\magnification=1200
\font\title=cmbx10 scaled\magstep2
\font\sc=cmcsc10

\font\tenmsbm=msbm10
\font\eightmsbm=msbm8
\textfont 9=\tenmsbm
\scriptfont 9=\eightmsbm
\def\bb {\fam9 }
\def\Der{\mathop{\rm Der}}
\def\ad{\mathop{\rm ad}}
\def\char{\mathop{\rm char}}
\centerline{\title On the Restricted Lie Algebra Structure}\smallskip

\centerline{\title for the Witt Lie Algebra in Finite Characteristic}\medskip

\centerline{\sc Tyler EVANS and Dmitry FUCHS}
\centerline{Department of Mathematics}
\centerline{University of California}
\centerline{Davis, CA 95616 USA}\bigskip

{\bf 1. Introduction.} Let $p$ be a prime, $\bb F$ be a field of
characteristic
$p$, and $\overline{\bb F}$ be the algebraic closure of $\bb F$. Let 
$A={\bb F} [x]/(x^p-1)$. Notice that ${\bb F}[x]/(x^p-1)\cong{\bb 
F}[x]/(x^p)$, an
isomorphism may be established be the formula $x\leftrightarrow x-1$. Notice
also that $\displaystyle{d\over dx}(x^p-1)\subset(x^p-1)$ and $\displaystyle
{d\over dx}(x^p)\subset(x^p)$, hence, the operator $\displaystyle{d\over dx}$
is also defined in $A$; the isomorphism described above commutes with
$\displaystyle{d\over dx}$. Below, we abbreviate the notation $\displaystyle
{d\over dx}$ to $\partial$.\smallskip

The Lie algebra $W=\Der A$ is called the Witt algebra. It consists of
``vector fields'' $f\partial,\ f\in A$. In particular, $\dim_{\bb F}W=\dim_{\bb
F}A=p.$

As any Lie algebra of derivations of a commutative algebra over $\bb F$, $W$
has a canonical structure of a restricted Lie algebra. Recall that a
restricted Lie algebra is a Lie algebra over $\bb F$ with an additional unary
(in general, non-linear) operation $g\mapsto g^{[p]}$ satisfying the
conditions
$$(\lambda g)^{[p]}=\lambda^pg^{[p]}\, (\lambda\in{\bb F}),\, \ad(g^{[p]})=
(\ad g)^p,$$ $$(g+h)^{[p]}=g^{[p]}+h^{[p]}+\sum_{i=1}^{p-1}s_i(g,h)$$where
$s_i(g,h)$ is the coefficient of $\lambda^{i-1}$ in $(\ad(\lambda g+h))^{p-1}
(h)$ times $i^{-1}$ modulo $p$; in particular, $[g^{[p]},g]=0$ for any $g$
(see details in [1]). In $\Der A$, $g^{[p]}=g^p$ (obviously, if $g\in\Der A$,
then $g^p=g\circ\dots\circ g\in\Der A$.)

Although the operation $g\mapsto g^{[p]}$ does not need to be linear, it is
fully determined by its values on any basis of the Lie algebra. In
particular,
in $W$,$$(x\partial)^{[p]}=x\partial,\, (x^k\partial)^{[p]}=0\ {\rm for}\
k=0,2,3,\dots,p-1$$(these formulas hold whether $A$ is regarded as
${\bb F}[x]/(x^p-1)$ or ${\bb F}[x]/(x^p)$; the same is true for Theorem 1
below). We will give, however, a more detailed description of the operation
$g\mapsto g^{[p]}$ in $W$.\smallskip

{\sc Theorem 1}. (a) {\it For any $f\in A$,
$(f\partial)^{[p]}=C(f)f\partial$,
where $C(f)$ is a constant (depending on $f$).}

(b) {\it For any $f\in A$, $\partial(f\partial(\dots(f\partial f)
\dots))$ with $p-1$ $\partial$'s (and $p-1$ $f$'s) is a constant, and this
constant is equal to $C(f)$}.

(c) {\it For any $f\in A$, $\partial^{p-1}(f^{p-1})$ is a constant, and
this constant is equal to $-C(f)$}.\smallskip

Since Parts (a) and (b) of Theorem 1 are very simple (see Section 2), it is
Part (c), or rather the equivalence between (b) and (c), that is the main
result of this paper. Moreover, this result is majorated by three (actually
equivalent) combinatorial theorems. We state these theorems in Section 3 and
prove them in Section 4.\smallskip

{\bf 2. Proof of Parts (a) and (b) of Theorem 1.} The Lie algebra $W$ has
rank one, in the sense that there exists a non-empty Zariski open subset
$U\subset W$ such that if $g\in U$ and $[h,g]=0$, then $h\in{\bb F}g$.  
Since $[g^{[p]},g]=0$, $g^{[p]}\in{\bb F}g$, at least for $g\in U$. But
since $g\mapsto g^{[p]}$ is an algebraic map, this implies that
$g^{[p]}\in{\bb F}g$ for all $g\in W$. (Strictly speaking, for this
purpose $\bb F$ should be infinite, but, if necessary, we can extend $\bb
F$ to $\overline{\bb F}$ and take $U$ in the extended $W$.) This proves
(a): $(f\partial)^{[p]}=C(f)f\partial$ for some algebraic function
$C\colon W\to{\bb F}$. To prove (b), apply the both sides of this equality
to $x\in A$. We have: $f\partial(f\partial(\dots(f\partial f)
\dots))=C(f)f$, which shows that $\partial(f\partial(\dots(f\partial f)
\dots))=C(f)$, at least if $f\in A$ is not a zero divisor. Hence, the last
equality holds for any $f$ in a non-empty Zariski open subset of $A$ (for
example, if $\tilde f\in{\bb F}[x]$ has a non-zero constant term, then the
image of $\tilde f$ under the projection ${\bb F}[x]\to{\bb F}[x]/(x^p)=A$
is not a zero divisor in $A$). Since the set of those $f\in A$ for which
our equality does not hold is also open, it should be empty: any two
non-empty Zariski open subsets of an affine space over an infinite field
overlap, and $\bar{\bb F}$, if not $\bb F$ itself, is infinite. Thus, the
equality holds for all $f\in A$, which is the statement of (b).\smallskip

{\bf 3. Three combinatorial theorems.} In this section, we state Theorems
2,3, and 4, and show that each of them implies Theorem 1 (c). In addition
to that, Theorems 2 and 3 are equivalent to each other and imply Theorem 
4.

First, we consider arbitrary finite words in a two-letter alphabet
$(\partial, f)$ ending with $f$. (We do not specify, what $f$ and
$\partial$ are; for example, $f$ may be a $C^\infty$ function in one
variable and $\partial$ the derivative.) This word may be regarded as an
integral linear combination of differential monomials $f^{(k_1)}\dots
f^{(k_m)}$. For example, $\partial f \partial\partial
f=(ff'')'=f'f''+ff'''$.\smallskip

{\sc Theorem 2.} {\it For any prime} $p$,$$(\partial f)^{p-1}\equiv-\partial
^{p-1}f^{p-1}\bmod p.$$

{\sc Example:}$$\eqalign{(\partial f)^4&=\partial f\partial f\partial
f\partial f=(f')^4+11f(f')^2f''+4f^2(f'')^2+7f^2f'f'''+f^3f^{(4)},\cr
\partial^4f^4&=24(f')^4+144f(f')^2f''+36f^2(f'')^2+48f^2f'f'''+4f^3f^{(4)}.
\cr}$$We see that $(\partial f)^4\equiv-\partial^4f^4\bmod5$, as stated in
Theorem.\smallskip

Theorem 2 (with Theorems 1(a), (b)) implies Theorem 1(c): take $f\in{\bb
F}[x], \partial=\displaystyle{d\over dx}$, and project the equality
$(\partial f)^ {p-1}=-\partial^{p-1}f^{p-1}$ (which holds if $\char{\bb
F}=p$) onto $A$.

Theorem 2 may be reformulated as a congruence of symmetric polynomials, in
the following way.\smallskip

{\sc Theorem 3.} {\it In ${\bb Z}[t_1,\dots,t_{p-1}]$,$$\sum_{\sigma\in
S_{p-1}}t_{\sigma(1)}(t_{\sigma(1)}+t_{\sigma(2)})\dots(t_{\sigma(1)}+\dots+
t_{\sigma(p-1)})\equiv(t_1+\dots+t_{p-1})^{p-1}\bmod p,$$where $p$, as usual,
is a prime.}\smallskip

(No minus sign is this congruence, it is not a misprint!)\smallskip

Obviously, if 
$$t_1(t_1+t_2)\dots(t_1+\dots+t_{p-1})=\sum n_{k_1\dots
k_{p-1}} t_1^{k_1}\dots t_{p-1}^{k_{p-1}},$$then$$\partial f_{p-1}\partial
f_{p-2}\dots\partial f_1=\sum n_{k_1\dots k_{p-1}}f_1^{(k_1)}\dots 
f_{p-1}^{(k_{p-1})}.$$

Similarly, if $$(t_1+\dots+t_{p-1})^{p-1}=\sum m_{k_1\dots k_{p-1}}
t_1^{k_1}\dots t_{p-1}^{k_{p-1}},$$ then $$\partial^{p-1}(f_1\dots
f_{p-1})= \sum m_{k_1\dots k_{p-1}}f_1^{(k_1)}\dots f_{p-1}^{(k_{p-1})}.$$
Hence, the congruence in Theorem 3 is equivalent to$$\sum_{\sigma\in
S_{p-1}}\partial f_{\sigma(p-1)}\partial f_{\sigma(p-2)}\dots\partial
f_{\sigma(1)}\equiv \partial^{p-1}(f_1\dots f_{p-1})\bmod p$$which
becomes, after substituting $f_1=\dots=f_{p-1}=f$, $$(p-1)!(\partial
f)^{p-1}\equiv\partial ^{p-1}f^{p-1}\bmod p$$
(the last two congruences are, actually, equivalent.) Since 
$(p-1)!\equiv-1\bmod p$, the last congruence is that of Theorem 2. Thus, 
Theorems 2 and 3 are equivalent.

Our last combinatorial theorem concerns a certain function on Young diagrams.
To avoid drawing, we use the term {\it Young diagram} for a finite sequence
$(j_1,\dots,j_m)$ of integers with $j_1\ge\dots\ge j_m>0$. The sequence may
be empty ($m=0$). For a Young diagram $J=(j_1,\dots,j_m)$, we put
$N(J)=j_1+\dots
+j_m,\, m(J)=m$, and $n_k(J)=\#\{s\mid j_s=k\}$. Define a function $d$ on
Young diagrams recursively: $d(\emptyset)=1$, and if $J=(j_1,\dots,j_m), 
N(J)=N$, and $d(K)$ has been already defined for all Young diagrams $K$ 
with $N(K)=N-1$,
then
$$d(J)=\sum_{s,j_s>j_{s+1}}(N-j_s+1)n_{j_s}(J)d(j_1,\dots,j_{s-1},j_s-1,
j_{s+1},\dots,j_m).$$
(Here we put $j_{m+1}=0$ and if $s=m$ and $j_s=1$, then
$j_s-1$ is zero, and we simply delete this zero.)\smallskip

{\sc Theorem 4}. {\it If $N(J)=p-1$ (where $p$ is prime), then}$$d(J)\equiv1
\bmod p.$$

{\sc Examples:}$$\eqalign{d(\emptyset)&=1;\cr d(1)&=1\cdot1\cdot d(\emptyset)
=1;\cr d(1,1)&=2\cdot2\cdot d(1)=4,\ d(2)=1\cdot1\cdot d(1)=1;\cr d(1,1,1)&=
3\cdot3\cdot d(1,1)=36,\ d(2,1)=2\cdot1\cdot d(1,1)+3\cdot1\cdot d(2)=11,\cr
&d(3)=1\cdot1\cdot d(2)=1;\cr d(1,1,1,1)&=4\cdot4\cdot d(1,1,1)=576,\
d(2,1,1)=3\cdot1\cdot d(1,1,1)+4\cdot2\cdot d(2,1)=196,\cr
&d(2,2)=3\cdot2\cdot
d(2,1)=66,\ d(3,1)=2\cdot1\cdot d(2,1)+4\cdot1\cdot d(3)=26,\cr &d(4)=1\cdot1
\cdot d(3)=1.\cr}$$

We see that if $N(J)=2$, then $d(J)=4,1\equiv1\bmod3$, and if $N(J)=4$, then
$d(J)=576,196,66,26,1\equiv1\bmod5$.

Theorem 4 is equivalent to Theorem 2 restricted to the case, when $f$ is a
monic polynomial of degree $p-1$ (this case of Theorem 2 is sufficient for
proving Theorem 1(c)).

Indeed, let $f(x)=(x-\alpha_1)\dots(x-\alpha_{p-1})$ (where $\alpha_1,\dots,
\alpha_{p-1}\in\overline{\bb F}$). We put $x-\alpha_i=u_i$; thus, $f=u_1\dots
u_{p-1}$ and $\partial u_i=1$. Let $n\le p-1$. Then $(\partial f)^n=
\partial f\partial f\dots\partial f$ is a symmetric polynomial in $u_1,\dots,
u_{p-1}$ of total degree $n(p-2)$ and of degree $\le p-1$ with respect to
each
variable $u_i$. Let $J=(j_1,\dots,j_m)$ be a Young diagram with $N(J)=n$.
Then
an obvious induction based on the equality $(\partial f)^n=\partial(u_1\dots
u_{p-1}(\partial f)^{n-1})$ shows that the coefficient at$$u_1^{n-j_1}\dots
u_m^{n-j_m}u_{m+1}^n\dots u_{p-1}^n$$ in the polynomial $(\partial f)^n$ is
$d(J)$.

On the other hand, the coefficient at the same monomial in the polynomial
$\partial^n(f^n)$ is$$\eqalign{{n!\over j_1!\dots j_m!}\prod_{i=1}^mn(n-1)
\dots(n-j_i+1)={n!\over j_1!\dots
j_m!}{n!\over(n-j_1)!}&\dots{n!\over(n-j_m)!}
\cr &=n!{n\choose j_1}\dots{n\choose j_m}.\cr}$$Since $(p-1)!\equiv-1\bmod p$
and $\displaystyle{p-1\choose j}\equiv(-1)^j\bmod p$, the last quantity for
$n=p-1$ is$$(p-1)!{p-1\choose j_1}\dots{p-1\choose j_m}\equiv(-1)\cdot
(-1)^{j_1}\cdot\dots\cdot(-1)^{j_m}\bmod p,$$and $(-1)\cdot
(-1)^{j_1}\cdot\dots\cdot(-1)^{j_m}=(-1)^p=-1$ (if $p$ is odd; if $p=2$, then
$-1\equiv1\bmod p$). Thus, Theorem 2 with $f=u_i\dots u_{p-1}$ is equivalent
to
Theorem 4.

We conclude this section with three remarks concerning Theorem 4. First, we
will never mention this theorem again; certainly, it follows from the other
theorems of this section, but we do not have any direct proof for it. Still
we
think that it deserves  to be stated as one of the results of this paper.
Second, this Theorem may have some meaning in the representation theory of
symmetric groups, but this meaning evades us. Third, it is not hard to deduce
from Theorem 4 that the congruence $d(J)\equiv1\bmod p$ holds also if $N(J)=
p-2$ (check this for $p=3$ and 5 using the example after the statement of
Theorem 4). We leave this to the reader as an exercise.\smallskip

{\bf 4. Proofs.} We will prove Theorem 3 (using its relations to propositions
similar to Theorems 1 and 2). As we know, this will imply all the other
theorems of this paper.

Let $\widetilde W=\Der{\bb F}[x]$. This is an infinite dimensional restricted
Lie algebra. Elements of $\widetilde W$ are ``vector fields'' $f\partial,\,
f\in{\bb F}[x]$. The $p$-th power of a derivation $f\partial$ is also a
derivation: $(f\partial)^p=F\partial,\, F\in{\bb F}[x]$. Raising $f\partial$
to the power $p$, we get$$F_1\partial+F_2\partial^2+\dots+F_p\partial^p=
F\partial$$(where $F_1=f\cdot(\partial(f\partial(\dots(f\partial f)\dots))),
F_p=f^p$). Applying both sides of this equality to $(x-a)^k$, where $1<k<p$
and $a\in\overline{\bb F}$, and then setting $x=a$, we get$$F_k(a)\cdot k!=
0,$$ which shows that $F_2=\dots=F_{p-1}=0$. Since $\partial^p=0$ on ${\bb
F}[x]$, we see that$$F=F_1=fg,\ g=(\partial f)^{p-1}=\partial f\partial
f\dots
\partial f.$$But $[(f\partial)^{[p]},f\partial]=0$; hence, $[fg\partial,f
\partial]=(fgf'-ff'g-f^2g')\partial=-f^2g'\partial=0$, that is, $g'=0$ (for
$f
\ne0$, and therefore for any $f$). (Actually, this means that $g$ is a
polynomial in $x^p$, but we will not need this.)

Consider differential expression$$g(f)=\partial\partial f\partial f\dots
\partial f\ (p\ \partial{\rm 's},\ p-1\ f{\rm 's}).$$Polarize the restriction
of the form $f\mapsto(g(f))(a),\, a\in\overline{\bb F}$ of degree $p-1$ to
the
vector space of polynomials of degree $<p$. We get a symmetric $(p-1)$-linear
form 
$$G(f_1,\dots,f_{p-1})=\sum_{\sigma\in S_{p-1}}(\partial\partial f_
{\sigma(1)}\partial f_{\sigma(2)}\dots\partial f_{\sigma(p-1)})(a)$$which
is equal to 0, since $g'(f)=0$. As a differential expression, the right hand
side of the last equality is a linear combination of monomials $f_1^{(j_1)}
f_2^{(j_2)}\dots f_{p-1}^{(j_{p-1})}(a)$ with $j_1+j_2+\dots+j_{p-1}=p$, but
plugging $f_1=(x-a)^{i_1},\dots,f_{p-1}=(x-a)^{i_{p-1}}$ with all $i_1,\dots,
i_{p-1}$ between 0 and $p-1$ and equating the results to 0, we see that all
the monomials with $j_1<p,\dots,j_{p-1}<p$ have zero coefficients in $\bb F$
(that is, they are 0 modulo $p$). Since the coefficient at $f_1\dots f_{s-1}
f_s^{(p)}f_{s+1}\dots f_{p-1}$ is obviously $(p-2)!\equiv1\bmod p$, we arrive
at the conclusion:$$\sum_{\sigma\in S_{p-1}}\partial\partial f_{\sigma(1)}
\dots\partial f_{\sigma(p-1)}\equiv\sum_{s=1}^{p-1}f_1\dots f_{s-1}f_s^{(p)}
f_{s+1}\dots f_{p-1}\bmod p,$$which may be rewritten as$$\sum_{\sigma\in
S_{p-1}}t_{\sigma(1)}(t_{\sigma(1)}+t_{\sigma(2)})\dots(t_{\sigma(1)}+\dots+
t_{\sigma(p-1)})(t_1+\dots+t_{p-1})\equiv t_1^p+\dots+t_{p-1}^p\bmod p$$(the
last factor in the left hand side of the last formula arises from the first
$\partial$ in the left hand side of the previous formula). But $t_1^p+\dots+
t_{p-1}^p\equiv(t_1+\dots+t_{p-1})^p\bmod p$. Canceling $t_1+\dots+t_{p-1}$,
we obtain the congruence of Theorem 3.\medskip

\centerline{\sc References}\medskip

[1] {\sc Jacobson N.} {\it Lie Algebras.} John Wiley, NY, 1962.\bye